\documentclass[12pt]{article}
\usepackage{amssymb,amsmath}
\usepackage{amsthm}
\usepackage{cases}
\usepackage{amsfonts}
\usepackage{graphicx}
\usepackage{cite,color,xcolor}
\usepackage[left=2.6cm,right=2.6cm,top=2.3cm,bottom=2.3cm]{geometry}
\usepackage[colorlinks,citecolor=blue,urlcolor=blue]{hyperref}
\usepackage[utf8]{inputenc}

\newtheorem{theorem}{Theorem}[section]

\newtheorem{lemma}{Lemma}[section]

\newtheorem{remark}{Remark}[section]

\usepackage{txfonts}
\newcommand{\bal}{\begin{align}}
\newcommand{\bbal}{\begin{align*}}
\newcommand{\beq}{\begin{equation}}
\newcommand{\eeq}{\end{equation}}
\newcommand{\bca}{\begin{cases}}
\newcommand{\eca}{\end{cases}}

\newcommand{\pa}{\partial}
\newcommand{\fr}{\frac}

\newcommand{\dd}{\mathrm{d}}

\newcommand{\R}{\mathbb{R}}

\newcommand{\f}{\left}
\newcommand{\g}{\right}
\begin{document}
\bibliographystyle{plain}
\title{Zero-filter limit for the Camassa-Holm equation in Sobolev spaces}

\author{Jinlu Li$^{1}$, Yanghai Yu$^{2,}$\footnote{E-mail: lijinlu@gnnu.edu.cn; yuyanghai214@sina.com(Corresponding author); mathzwp2010@163.com} and Weipeng Zhu$^{3}$\\
\small $^1$ School of Mathematics and Computer Sciences, Gannan Normal University, Ganzhou 341000, China\\
\small $^2$ School of Mathematics and Statistics, Anhui Normal University, Wuhu 241002, China\\
\small $^3$ School of Mathematics and Big Data, Foshan University, Foshan, Guangdong 528000, China}

\date{\today}
\maketitle\noindent{\hrulefill}

{\bf Abstract:} The aim of this paper is to answer the question left in \cite{GL} (Math. Z. (2015) 281). We prove that the zero-filter limit of the Camassa-Holm equation is the Burgers equation in the same topology of Sobolev spaces as the initial data.



\vskip0mm\noindent{\hrulefill}

\section{Introduction}
The Camassa-Holm equations with fractional dissipation reads as follows
\begin{equation}\label{fC}
\begin{cases}
\pa_tm+2m\pa_xu+u\pa_xm+\nu\Lambda^{\gamma}m=0,\\
m=(1-\alpha^2\pa^2_x)u,\\
u(0,x)=u_0(x),
\end{cases}
\end{equation}
where  the constant $\alpha>0$ is a filter parameter, the
constants $\nu\geq0$ and $\gamma\in[0,2]$. The fractional power operator $\Lambda^{\gamma}$ is defined by Fourier multiplier with the symbol $|\xi|^{\gamma}$, namely,
  $\Lambda^{\gamma}u(x)=\mathcal{F}^{-1}\big(|\xi|^{\gamma}\mathcal{F}u(\xi)\big).
$
When $\nu=0$, \eqref{fC} becomes the classical Camassa-Holm (CH) equation
\begin{equation}\label{C}
\begin{cases}
\pa_tm+2m\pa_xu+u\pa_xm=0,\\
m=(1-\alpha^2\pa^2_x)u,\\
u(0,x)=u_0(x).
\end{cases}
\end{equation}
The CH equation was firstly proposed in the context of hereditary symmetries studied in \cite{Fokas} and then was derived explicitly as a water wave equation by Camassa-Holm \cite{Camassa}. Many aspects of the mathematical beauty of the CH equation have been exposed over the last two decades. Particularly, \eqref{C} is completely integrable \cite{Camassa,Constantin-P} with a bi-Hamiltonian structure \cite{Constantin-E,Fokas} and infinitely many conservation laws \cite{Camassa,Fokas}. Also, it admits exact peaked
soliton solutions (peakons) of the form $u(x,t)=ce^{-|x-ct|}$ with $c>0,$ which are orbitally stable \cite{Constantin.Strauss}. Another remarkable feature of the CH equation is the wave breaking phenomena: the solution remains bounded while its slope becomes unbounded in finite time \cite{Constantin,Escher2,Escher3}. It is worth mentioning that the peaked solitons present the characteristic for the travelling water waves of greatest height and largest amplitude and arise as solutions to the free-boundary problem for incompressible Euler equations over a flat bed, see Refs. \cite{Constantin-I,Escher4,Escher5} for the details. Because of the interesting and remarkable features as mentioned above, the Camassa-Holm equation has attracted much attention as a class of integrable shallow water wave equations in recent twenty years. Its systematic mathematical study was initiated in a series of papers by Constantin and Escher, see \cite{Escher1,Escher2,Escher3,Escher4,Escher5}.
When $\nu > 0$ and $\gamma= 2$, Xin-Zhang \cite{xin} proved that the classical viscous Camassa-Holm equation is globally well-posed.

Note that
\bal\label{g}
\quad\left(1-\alpha^2 \partial_x^2\right)^{-1} f=g*f,\quad \forall \;f \in L^2(\mathbb{R}),
\end{align}
where $g(x):=\frac{1}{2 \alpha} e^{-\frac{|x|}{\alpha}},\; x \in \mathbb{R}$ and $*$ denotes convolution, then $u=g * m$. Using this identity and applying the pseudodifferential operator $\left(1-\alpha^2 \partial_x^2\right)^{-1}$ to Eq. \eqref{fC}, one can rewrite Eq. \eqref{fC} as a quasi-linear nonlocal evolution equation of hyperbolic type, namely
\begin{align}\label{ch}
\begin{cases}
u_t+u \partial_xu+\nu\Lambda^{\gamma}u=-\partial_x \left(1-\alpha^2 \partial_x^2\right)^{-1}\left(u^2+\frac{\alpha^2}{2} (\partial_xu)^2\right),  &\quad (t,x)\in \R^+\times\R,\\
u(0,x)=u_0(x).
\end{cases}
\end{align}
When the filter parameter $\alpha=0$, Eq. \eqref{ch} becomes the Burgers equation
\begin{align}\label{bu}
\begin{cases}
u_t+3u \partial_xu+\nu\Lambda^{\gamma}u=0,  &\quad (t,x)\in \R^+\times\R,\\
u(0,x)=u_0(x).
\end{cases}
\end{align}
The Burgers equation \eqref{fb} with $\gamma=0$ and $\gamma=2$ has received an extensive amount of attention since the studies by Burgers in the 1940s. If $\gamma=0$, the equation is perhaps the most basic example of a PDE evolution leading to shocks. If $\gamma=2$, it provides an accessible model for studying the interaction between nonlinear and dissipative phenomena. Kiselev et al. \cite{Kiselev} gave a complete study for general $\gamma\in[0,2]$ for the periodic case. In particular, for the case $\gamma=1$, they proved the global well-posedness of the equation in the critical Hilbert space $H^{\fr12}(\mathbb{T})$ by using the method of modulus of continuity. Subsequently, Miao-Wu \cite{miao2009} proved the global well-posedness of the critical Burgers equation in critical Besov spaces $B^{1/p}_{p,1}(\R)$ with $p\in[1,\infty)$ with the help of Fourier localization technique and the method of modulus of continuity. For more results on the fractional Burgers equation and dispersive perturbations of Burgers equations, we refer the readers to see \cite{Alibaud,Dong,Karch,Linares,Molinet} and the references therein.

Formally, as $\alpha\to0$, the solution of the Camassa-Holm equation \eqref{ch} converges to the
solution of the Burgers equation \eqref{bu}. More precisely, Gui-Liu \cite{GL} proved that the
solutions of the dissipative Camassa-Holm equation ($\nu>0, \gamma\in(0,1]$) does converge, at least locally, to the
one of the dissipative Burgers equation as the filter parameter $\alpha$ tends to zero in the lower regularity Sobolev spaces. {\it However, the question about whether the zero-filter limit of solutions of the classical
Camassa-Holm equation is a solution to the inviscid Burgers equation, as mentioned in \cite{GL} (see p.997), is still an open problem.}

In this paper, we consider the problem of the zero-filter limit $(\alpha\to0)$  of the following Camassa-Holm equation
\begin{align}\label{fb}
\begin{cases}
u_t+u \partial_xu=-\partial_x \left(1-\alpha^2 \partial_x^2\right)^{-1}\left(u^2+\frac{\alpha^2}{2} (\partial_xu)^2\right), \\
u(0,x)=u_0(x)\in H^s(\R).
\end{cases}
\end{align}
Equivalently,
\begin{align}\label{fb1}
\begin{cases}
u_t+3u \partial_xu=-\alpha^2\partial^3_x \left(1-\alpha^2 \partial_x^2\right)^{-1}(u^2)-\frac{\alpha^2}{2}\partial_x \left(1-\alpha^2 \partial_x^2\right)^{-1} (\partial_xu)^2,\\
u(0,x)=u_0(x)\in H^s(\R).
\end{cases}
\end{align}
In this paper, we shall address the above
open question. Precisely speaking, we shall prove that the solution of \eqref{fb} converges to the solution of the following inviscid Burgers equation in the topology of Sobolev spaces
\begin{align}\label{ib}
\begin{cases}
u_t+3u \partial_xu=0, &\quad (t,x)\in \R^+\times\R,\\
u(0,x)=u_0(x)\in H^s(\R).
\end{cases}
\end{align}
Our main result is the following:
\begin{theorem}\label{th1} Let $s>\frac32$ and $\alpha \in(0,1)$. Assume that the initial data $u_0\in H^s(\mathbb{R})$. Let $\mathbf{S}_{t}^{\mathbf{\alpha}}(u_0)$ and $\mathbf{S}_{t}^{0}(u_0)$ be the smooth solutions of \eqref{fb} and \eqref{ib} with the initial data $u_0$ respectively. Then there exists a time $T=T(\|u_0\|_{H^s})>0$ such that  $\mathbf{S}_{t}^{\mathbf{\alpha}}(u_0),\mathbf{S}_{t}^{0}(u_0)\in \mathcal{C}([0,T];H^s)$ and
$$
\lim_{\alpha\to 0}\left\|\mathbf{S}_{t}^{\mathbf{\alpha}}(u_0)-\mathbf{S}_{t}^{0}(u_0)\right\|_{L^\infty_TH^{s}}=0.
$$
\end{theorem}
\begin{remark}\label{re1}
Compared with the result ($\nu>0$) in the weak topology Sobolev spaces given by Gui-Liu in \cite{GL}, our Theorem \ref{th1} holds for $\nu=0$ (zero dissipative) and seems to be optimal in the sense of that the  convergence space is  the solution spaces of both the Camassa-Holm equation and inviscid Burgers equation. We also would like to emphasize that, Theorem \ref{th1} holds for any $\nu\geq0$.
\end{remark}
\begin{remark}\label{re2}
Motivated by the Bona–Smith method \cite{Bona} (see also \cite{GLY}), we decompose the difference of $\mathbf{S}_{t}^{\mathbf{\alpha}}(u_0)$ and $\mathbf{S}_{t}^{0}(u_0)$ as follows:
\bal\label{vw}
\mathbf{S}_{t}^{\mathbf{\alpha}}(u_0)-\mathbf{S}_{t}^{0}(u_0)&=\mathbf{S}_{t}^{\mathbf{\alpha}}(u_0)-\mathbf{S}_{t}^{\alpha}(S_nu_0)
+\mathbf{S}_{t}^{\alpha}(S_nu_0)-\mathbf{S}_{t}^{0}(S_nu_0)
+\mathbf{S}_{t}^{0}(S_nu_0)-\mathbf{S}_{t}^{0}(u_0),
\end{align}
where $S_n$ denotes the frequency localization operator defined in Lemma \ref{c}.
The key idea is to show that
$$\|\mathbf{S}_{t}^{\mathbf{\alpha}}(u_0)-\mathbf{S}_{t}^{\alpha}(S_nu_0)\|_{H^{s}}\leq C\|(\mathrm{Id}-S_n)u_0\|_{H^{s}}, \quad \forall \alpha\in[0,1)$$
and
\begin{align*}
\|\mathbf{S}_{t}^{\alpha}(S_nu_0)-\mathbf{S}_{t}^{0}(S_nu_0)\|_{H^{s}} \leq C \alpha 2^{\frac32n}.
\end{align*}
For more details see {\bf Step 2} and {\bf Step 3} in Section \ref{sec3}.
\end{remark}
\begin{remark}\label{re3}
We should mention that, by the idea in Remark \ref{re2}, our Theorem \ref{th1} holds for the Besov spaces $B^s_{p,r}(\R)$ with $s>\max\{\fr32,1+\frac1p\}$ and $(p,r)\in (1,\infty)\times [1,\infty)$. In order to elucidate the main idea, we do not pursue the general case in the current paper and leave it to the interested readers.
\end{remark}
\section{Preliminaries}\label{sec2}
{\bf Notation}\; Throughout this paper, we will denote by $C$ any positive constant independent of the parameter $\alpha$, which may change from line to line.
Given a Banach space $X$, we denote its norm by $\|\cdot\|_{X}$. For $I\subset\R$, we denote by $\mathcal{C}(I;X)$ the set of continuous functions on $I$ with values in $X$. Sometimes we will denote $L^p(0,T;X)$ by $L_T^pX$.
For all $f\in \mathcal{S}'$, the Fourier transform  $\widehat{f}$ is defined by
$$
\widehat{f}(\xi)=\int_{\R}e^{-ix\xi}f(x)\dd x \quad\text{for any}\; \xi\in\R.
$$
For $s\in\R$, we denote the operator $J^{s}:=(1-\pa_x^2)^{\fr{s}2}$ which is defined by
$$\widehat{J^s f}(\xi)=(1+|\xi|^2)^{\fr{s}2}\widehat{f}(\xi).$$
For $s\in\R$, the nonhomogeneous Sobolev space is defined by
$$\|f\|^2_{H^s}=\|J^{s}f\|^2_{L^2}=\int_{\R}(1+|\xi|^2)^s|\widehat{f}(\xi)|^2\dd \xi.$$
Next, we introduce some known results for later proof.
\begin{lemma}[\cite{B}]\label{le2}
For $s>0$, $H^s(\R)\cap L^\infty(\R)$ is an algebra.
Moreover, we have for any $u,v \in H^s(\R)\cap L^\infty(\R)$
\begin{align*}
&\|uv\|_{H^s(\R)}\leq C\big(\|u\|_{H^s(\R)}\|v\|_{L^\infty(\R)}+\|v\|_{H^s(\R)}\|u\|_{L^\infty(\R)}\big).
\end{align*}
In particular, for $s>\frac12$, due to the fact $H^s(\R)\hookrightarrow L^\infty(\R)$, then we have
\begin{align*}
&\|uv\|_{H^s(\R)}\leq C\|u\|_{H^s(\R)}\|v\|_{H^s(\R)}.
\end{align*}
\end{lemma}
\begin{lemma}[\cite{Kato}]\label{le3}
Let $s>0$ and $f,g\in {\rm{Lip}}\cap H^s(\R)$ and $g\in L^\infty\cap H^{s-1}(\R)$. Then we have
\bbal
\|[J^s,f]g\|_{L^2}\leq C\big(\|\partial_x f\|_{L^\infty}\|g\|_{H^{s-1}}+\|f\|_{H^s}\|g\|_{L^\infty}\big).
\end{align*}
\end{lemma}
\begin{lemma}[\cite{B}]\label{c}
Let the inhomogeneous low-frequency cut-off operator $S_{n}$ is defined by
$
S_n u:=\sum\limits_{-1\leq q\leq n-1}{\Delta}_qu.
$
For any $u \in H^s(\R)$ with $s>\fr32$, we have
\bbal
&\lim_{n\to+\infty}\|S_nu-u\|_{H^{s}(\R)}=0.
\end{align*}
\end{lemma}
\section{Proof of Theorem \ref{th1}}\label{sec3}
 \quad We divide the proof of Theorem \ref{th1} into three steps.

 {\bf Step 1: Uniform bound w.r.s $\alpha\in(0,1)$ of $\mathbf{S}_{t}^{\mathbf{\alpha}}(u_0)$ in $H^s$}.

For fixed $\alpha>0$, by the classical local well-posedness result, we known that there exists a $T_\alpha=T(\|u_0\|_{H^s},\alpha)>0$ such that the Camassa-Holm has a unique solution $\mathbf{S}_{t}^{\mathbf{\alpha}}(u_0)\in\mathcal{C}([0,T_\alpha];H^s)$.

We shall prove that $\exists\; T=T(\|u_0\|_{H^s})>0$ such that $T\leq T_{\alpha}$ and there exists $C>0$ independent of $\alpha$ such that
\begin{align}\label{m1}
\|\mathbf{S}_{t}^{\mathbf{\alpha}}(u_0)\|_{L_T^{\infty} H^s} \leq C, \quad \forall \alpha \in[0,1).
\end{align}
Moreover, if $u_0 \in H^\gamma$ for some $\gamma>s$, then there exists $C_2>0$ independent of $\alpha$ such that
\begin{align}\label{m2}
\|\mathbf{S}_{t}^{\mathbf{\alpha}}(u_0)\|_{L_T^{\infty} H^\gamma} \leq C_2\left\|u_0\right\|_{H^\gamma} .
\end{align}
We just prove \eqref{m1} since \eqref{m2} is similar. To simplify notation, we set $u=\mathbf{S}_{t}^{\mathbf{\alpha}}(u_0)$.
Applying the operator $J^{s}uJ^{s}$ to $\eqref{ch}$ and integrating the resulting over $\R$, we obtain
\begin{align}
\frac12\frac{\dd }{\dd t}\|u\|^2_{{H}^{s}}&=\int_{\R}\pa_xu|J^{s}u|^2\dd x-\int_{\R}[J^s,u]\pa_xu\cdot J^su\dd x\label{y1}\\
&\quad-2\int_{\R}(1-\alpha^2\pa^2_x)^{-1}J^s(uu_x)\cdot J^su\dd x\label{y2}\\
&\quad
-\frac{\alpha^2}{2}\int_{\R}\pa_x(1-\alpha^2\pa^2_x)^{-1}J^s(\partial_xu)^2\cdot J^su\dd x.\label{y3}
\end{align}
To bound \eqref{y1}, by the classical commutator estimation (see Lemma \ref{le3}), it is easy to obtain
\bbal
|\eqref{y1}|
&\leq C\f(\|\pa_xu\|_{L^\infty}\|u\|^2_{H^s}+\|[J^s,u]\pa_xu\|_{L^2}\|u\|_{H^s}\g)\leq C\|\pa_xu\|_{L^\infty}\|u\|^2_{H^s},
\end{align*}
To bound \eqref{y2}, notice that
\begin{align}
\int_{\R}(1-\alpha^2\pa^2_x)^{-1}J^s(uu_x)J^su\dd x&=\int_{\R}(1-\alpha^2\pa^2_x)^{-1}( u\pa_xJ^su)\cdot J^su\dd x
\label{y4}\\& \quad +\int_{\R}(1-\alpha^2\pa^2_x)^{-1}[J^s,u]\pa_xu\cdot J^su\dd x,\label{y5}
\end{align}
and letting $v=(1-\alpha^2\pa^2_x)^{-1}J^su$, then
\bbal
\eqref{y4}&=\int_{\R} u(1-\alpha^2\pa^2_x)\pa_xv\cdot v\dd x
\\&=-\frac12\int_{\R} \pa_xu\cdot v^2\dd x+\int_{\R} \pa_xu(\alpha^2\pa^2_xv)\cdot v\dd x-\fr12\int_{\R} \pa_xu\cdot (\alpha\pa_xv)^2\dd x\\
&\leq C\|\pa_xu\|_{L^\infty}\f(\|v\|^2_{L^2}+\|\alpha^2\pa^2_xv\|^2_{L^2}+\|\alpha\pa_xv\|^2_{L^2}\g)\\
&\leq C\|\pa_xu\|_{L^\infty}\|u\|^2_{H^s},
\end{align*}
where we have used \eqref{g} and the convolution inequality.
Also,
\bbal
\eqref{y5}&=\int_{\R} [J^s,u]\pa_xu\cdot v\dd x
\leq \|[J^s,u]\pa_xu\|_{L^2}\|v\|_{L^2}
\leq C\|\pa_xu\|_{L^\infty}\|u\|^2_{H^s},
\end{align*}
which implies that
\bbal
|\eqref{y2}|\leq C\|\pa_xu\|_{L^\infty}\|u\|^2_{H^s}.
\end{align*}
To bound \eqref{y3}, notice that
\bbal
\eqref{y3}&=\frac{\alpha^2}{2}\int_{\R}J^{s-1}(\partial_xu)^2\cdot\pa_xJ v\dd x\\
&\leq \alpha^2\|(\partial_xu)^2\|_{H^{s-1}}\|\pa_xJ v\|_{L^2}\\
&\leq C\|\pa_xu\|_{L^\infty}\|u\|^2_{H^s}.
\end{align*}
Combining the above yields that
\bbal
\frac{\dd}{\dd t}\|\mathbf{S}_{t}^{\mathbf{\alpha}}(u_0)\|^2_{H^s}\leq C\|\pa_x\mathbf{S}_{t}^{\mathbf{\alpha}}(u_0)\|_{L^{\infty}}\|\mathbf{S}_{t}^{\mathbf{\alpha}}(u_0)\|^2_{H^s}\leq C\|\mathbf{S}_{t}^{\mathbf{\alpha}}(u_0)\|^3_{H^s}.
\end{align*}
Thus, by continuity arguments there exists a time $T=T(\|u_0\|_{H^s})>0$ such that \eqref{m1} holds uniformly w.r.s $\alpha\in(0,1)$.

{\bf Step 2: Estimations of $\|\mathbf{S}_{t}^{\mathbf{\alpha}}(u_0)-\mathbf{S}_{t}^{\alpha}(S_nu_0)\|_{H^s}$ and $\|\mathbf{S}_{t}^{0}(u_0)-\mathbf{S}_{t}^{0}(S_nu_0)\|_{H^s}$}.

 Denoting $$\mathbf{v}(t)=\mathbf{S}_{t}^{\mathbf{\alpha}}(u_0)-\mathbf{S}_{t}^{\alpha}(S_nu_0)\quad\text{and}\quad \mathbf{v}|_{t=0}=(\mathrm{Id}-S_n)u_0,$$ we infer that $\mathbf{v}$ satisfies
$$
\partial_{t} \mathbf{v}+\mathbf{S}_{t}^{\mathbf{\alpha}}(u_0)\partial_{x} \mathbf{v}=-\mathbf{v}\partial_{x} \mathbf{S}_{t}^{\alpha}(S_nu_0) +\mathcal{B}(\mathbf{v}, \mathbf{S}_{t}^{\mathbf{\alpha}}(u_0)+\mathbf{S}_{t}^{\alpha}(S_nu_0)),
$$
where
$$
\mathcal{B}:(f, g) \mapsto \partial_x\left(1-\alpha^2 \partial_x^2\right)^{-1}\left(f g+\frac{\alpha^2}{2} \partial_{x} f \partial_{x} g\right).
$$
Notice that
\bal
&\quad \ \int_{\R}J^s\f(\mathcal{B}(\mathbf{v}, \mathbf{S}_{t}^{\mathbf{\alpha}}(u_0)+\mathbf{S}_{t}^{\alpha}(S_nu_0))\g)\cdot J^s\mathbf{v}\dd x \nonumber
\\&=\frac12\int_{\R}J^s\f(\alpha^2\pa_x(1-\alpha^2\pa^2_x)^{-1}(\mathbf{v}_x \pa_x[\mathbf{S}_{t}^{\mathbf{\alpha}}(u_0)+\mathbf{S}_{t}^{\alpha}(S_nu_0)])\g)\cdot J^s\mathbf{v}\dd x \label{lyz+1}
\\&\quad +2\int_{\R}(1-\alpha^2\pa^2_x)^{-1} J^s(\pa_x\mathbf{S}_{t}^{\alpha}(S_nu_0)\mathbf{v})\cdot J^s\mathbf{v}\dd x\label{lyz+2}\\
&\quad +2\int_{\R}(1-\alpha^2\pa^2_x)^{-1}J^s\f[\f(\mathbf{v}+\mathbf{S}_{t}^{\alpha}(S_nu_0)\g)\mathbf{v}_x)\g]\cdot J^s\mathbf{v}\dd x. \label{lyz+3}
\end{align}
To bound \eqref{lyz+1} and \eqref{lyz+2}, we easily obtain
\bal\label{100}
|\eqref{lyz+1}|+|\eqref{lyz+2}|&\leq C\f(\|\mathbf{S}_{t}^{\mathbf{\alpha}}(u_0)\|_{H^s}
+\|\mathbf{S}_{t}^{\alpha}(S_nu_0)\|_{H^s}\g)\|\mathbf{v}\|^2_{H^s}
+C\|\pa_x\mathbf{S}_{t}^{\alpha}(S_nu_0)\|_{H^{s}}\|\mathbf{v}\|_{H^{s-1}}\|\mathbf{v}\|_{H^{s}}.
\end{align}
Similar argument as \eqref{y2}, we have
\bal\label{101}
|\eqref{lyz+3}|\leq C\f(\|\mathbf{v}\|_{H^s}+\|\mathbf{S}_{t}^{\alpha}(S_nu_0)\|_{H^s}\g)\|\mathbf{v}\|^2_{H^s}\leq C\f(\|\mathbf{S}_{t}^{\alpha}(u_0)\|_{H^s}+\|\mathbf{S}_{t}^{\alpha}(S_nu_0)\|_{H^s}\g)\|\mathbf{v}\|^2_{H^s}.
\end{align}
Following the same procedure as that in Step 1, we deduce from \eqref{100}-\eqref{101} that
\begin{align}\label{l0}
\frac12\frac{\dd }{\dd t}\|\mathbf{v}\|^2_{{H}^{s}}&=\int_{\R}\pa_x\mathbf{S}_{t}^{\mathbf{\alpha}}(u_0)|J^{s}\mathbf{v}|^2\dd x-\int_{\R}[J^s,\mathbf{S}_{t}^{\mathbf{\alpha}}(u_0)]\pa_x\mathbf{v}\cdot J^s\mathbf{v}\dd x\nonumber\\
&\quad-\int_{\R}J^s(\mathbf{v}\partial_{x} \mathbf{S}_{t}^{\alpha}(S_nu_0) )\cdot J^s\mathbf{v}\dd x\nonumber\\
&\quad
+\int_{\R}J^s\f(\mathcal{B}(\mathbf{v}, \mathbf{S}_{t}^{\mathbf{\alpha}}(u_0)+\mathbf{S}_{t}^{\alpha}(S_nu_0))\g)\cdot J^s\mathbf{v}\dd x\nonumber\\
&\leq C\f(\|\pa_x\mathbf{S}_{t}^{\mathbf{\alpha}}(u_0)\|_{L^\infty}+\|\mathbf{S}_{t}^{\mathbf{\alpha}}(u_0)\|_{H^s}+\|\mathbf{S}_{t}^{\alpha}(S_nu_0)\|_{H^s}\g)\|\mathbf{v}\|^2_{H^s}\nonumber\\
&\quad+ C\|\pa_x\mathbf{S}_{t}^{\alpha}(S_nu_0)\|_{H^{s}}\|\mathbf{v}\|_{H^{s-1}}\|\mathbf{v}\|_{H^{s}}.
\end{align}
Using \eqref{m1} and \eqref{m2}, the above \eqref{l0} reduces to
\begin{align}\label{l3}
\frac{\dd }{\dd t}\|\mathbf{v}\|_{{H}^{s}}
&\leq C\|\mathbf{v}\|_{H^s}+ C2^n\|\mathbf{v}\|_{H^{s-1}}.
\end{align}
To close \eqref{l3}, we have to estimate $\|\mathbf{v}\|_{H^{s-1}}$.
\begin{align*}
\frac12\frac{\dd }{\dd t}\|\mathbf{v}\|^2_{{H}^{s-1}}&=\int_{\R}\pa_x\mathbf{S}_{t}^{\mathbf{\alpha}}(u_0)|J^{s-1}\mathbf{v}|^2\dd x-\int_{\R}[J^{s-1},\mathbf{S}_{t}^{\mathbf{\alpha}}(u_0)]\pa_x\mathbf{v}\cdot J^{s-1}\mathbf{v}\dd x\nonumber\\
&\quad-\int_{\R}J^{s-1}(\mathbf{v}\partial_{x} \mathbf{S}_{t}^{\alpha}(S_nu_0) )\cdot J^{s-1}\mathbf{v}\dd x\nonumber\\
&\quad
+\int_{\R}J^{s-1}\f(\mathcal{B}(\mathbf{v}, \mathbf{S}_{t}^{\mathbf{\alpha}}(u_0)+\mathbf{S}_{t}^{\alpha}(S_nu_0))\g)\cdot J^{s-1}\mathbf{v}\dd x\nonumber\\
&\leq C\f(\|\pa_x\mathbf{S}_{t}^{\mathbf{\alpha}}(u_0)\|_{L^\infty}+\|\mathbf{S}_{t}^{\mathbf{\alpha}}(u_0)\|_{H^{s}}+\|\mathbf{S}_{t}^{\alpha}(S_nu_0)\|_{H^{s}}\g)\|\mathbf{v}\|^2_{H^{s-1}}.
\end{align*}
Applying Gronwall's inequality yields that for $t\in [0,T]$
\bal\label{l4}
\|\mathbf{v}(t)\|_{H^{s-1}}\leq C\|(\mathrm{Id}-S_n)u_0\|_{H^{s-1}}\leq C2^{-n}\|(\mathrm{Id}-S_n)u_0\|_{H^{s}}.
\end{align}
Inserting \eqref{l4} into \eqref{l3} and applying Gronwall's inequality, we obtain that for $t\in [0,T]$
\bbal
\|\mathbf{S}_{t}^{\mathbf{\alpha}}(u_0)-\mathbf{S}_{t}^{\alpha}(S_nu_0)\|_{H^{s}}\leq C\|(\mathrm{Id}-S_n)u_0\|_{H^{s}}, \quad \forall \alpha\in[0,1).
\end{align*}

{\bf Step 3: Estimation of $\|\mathbf{S}_{t}^{\alpha}(S_nu_0)-\mathbf{S}_{t}^{0}(S_nu_0)\|_{H^s}$}.

We can find that $\mathbf{S}_{t}^{\alpha}(S_nu_0)$ satisfies the following equation
\bbal
\pa_t\mathbf{S}_{t}^{\alpha}(S_nu_0)+3\mathbf{S}_{t}^{\alpha}(S_nu_0)\pa_x\mathbf{S}_{t}^{\alpha}(S_nu_0)
&=-\frac12\alpha^2\pa_x(1-\alpha^2\pa^2_x)^{-1}[\pa_x\mathbf{S}_{t}^{\alpha}(S_nu_0)]^2
\\&\quad -\alpha^2\pa^3_x(1-\alpha^2\pa^2_x)^{-1}[\mathbf{S}_{t}^{\alpha}(S_nu_0)]^2.
\end{align*}
Denoting $$\mathbf{w}(t)= \mathbf{S}_{t}^{\alpha}(S_nu_0)-\mathbf{S}_{t}^{0}(S_nu_0)\quad\text{and}\quad \mathbf{w}|_{t=0}=0,$$ we infer that $\mathbf{w}$ satisfies
$$
\left\{\begin{array}{l}
\partial_t\mathbf{w}+3\mathbf{S}_{t}^{0}(S_nu_0) \partial_x \mathbf{w}=-3\mathbf{w} \partial_x \mathbf{S}_{t}^{\alpha}(S_nu_0)-\mathbf{I}, \\
\mathbf{w}|_{t=0}=0,
\end{array}\right.
$$
where we denote
\bbal
\mathbf{I}:=\alpha^2\partial_x\left(1-\alpha^2 \partial_x^2\right)^{-1}\f[\partial^2_x \left([\mathbf{S}_{t}^{\alpha}(S_nu_0)]^2\right)+\frac{1}{2}\left[\pa_x\mathbf{S}_{t}^{\alpha}(S_nu_0)\right]^2\g].
\end{align*}
Taking the similar argument with \eqref{l0}, we have
\begin{align}
\frac12\frac{\dd }{\dd t}\|\mathbf{w}\|^2_{{H}^{s-1}}&=\int_{\R}\pa_x\mathbf{S}_{t}^{0}(S_nu_0)|J^{s-1}\mathbf{w}|^2\dd x-\int_{\R}[J^{s-1},\mathbf{S}_{t}^{0}(S_nu_0)]\pa_x\mathbf{w}\cdot J^{s-1}\mathbf{w}\dd x\label{w1}\\
&\quad-3\int_{\R}J^{s-1}\f(\mathbf{w} \partial_x \mathbf{S}_{t}^{\alpha}(S_nu_0)\g)\cdot J^{s-1}\mathbf{w}\dd x\label{w2}\\
&\quad-\int_{\R}J^{s-1}\mathbf{I}\cdot J^{s-1}\mathbf{w}\dd x\label{w3}.
\end{align}
Obviously,
\begin{align*}
|\eqref{w1}|+|\eqref{w2}|&\leq C(\|\mathbf{S}_{t}^{0}(S_nu_0)\|_{H^s}+\|\mathbf{S}_{t}^{\alpha}(S_nu_0)\|_{H^s})\|\mathbf{w}\|^2_{H^{s-1}},\\
|\eqref{w3}|&\leq C\|\mathbf{I}\|_{H^{s-1}}\|\mathbf{w}\|_{H^{s-1}}\\
&\leq C\alpha^2\|\mathbf{S}_{t}^{\alpha}(S_nu_0)\|_{W^{1,\infty}}\|\mathbf{S}_{t}^{\alpha}(S_nu_0)\|_{H^{s+2}}\|\mathbf{w}\|_{H^{s-1}}\\
&\leq C\alpha^22^{2n}\|\mathbf{w}\|_{H^{s-1}}.
\end{align*}
Gathering the above estimates,  we deduce that
\begin{align*}
\|\mathbf{w}(t)\|_{H^{s-1}}
&\leq C \int_{0}^{t}\|\mathbf{w}\|_{H^{s-1}} \dd\tau+C\alpha^22^{2n},
\end{align*}
which along with Gronwall's inequality implies
\bbal
\|\mathbf{w}(t)\|_{H^{s-1}}\leq C\alpha^22^{2n}.
\end{align*}
Then, we get for $t\in[0,T]$
\begin{align*}
\|\mathbf{w}(t)\|_{H^{s}} &\leq \|\mathbf{w}(t)\|^{\frac12}_{H^{s-1}}\|\mathbf{w}(t)\|^{\frac12}_{H^{s+1}} \leq C \alpha 2^{\frac32n}.
\end{align*}
Due to \eqref{vw} and using {\bf Step 1}-{\bf Step 3}, we have for $t\in[0,T]$
\bbal
&\|\mathbf{S}_{t}^{\mathbf{\alpha}}(u_0)-\mathbf{S}_{t}^{0}(u_0)\|_{H^{s}}\\
&\leq \|\mathbf{S}_{t}^{\mathbf{\alpha}}(u_0)-\mathbf{S}_{t}^{\alpha}(S_nu_0)\|_{H^{s}}+\|\mathbf{S}_{t}^{\alpha}(S_nu_0)-\mathbf{S}_{t}^{0}(S_nu_0)\|_{H^s}
+\|\mathbf{S}_{t}^{0}(S_nu_0)-\mathbf{S}_{t}^{0}(u_0)\|_{H^{s}}\\
&\leq C\|(\mathrm{Id}-S_n)u_0\|_{H^{s}}+C2^{\frac32n}\alpha.
\end{align*}
Using Lemma \ref{c} to the above inequality enables us to complete the proof of Theorem \ref{th1}.

\section*{Acknowledgments}
J. Li is supported by the National Natural Science Foundation of China (11801090 and 12161004) and Jiangxi Provincial Natural Science Foundation (20212BAB211004). Y. Yu is supported by the National Natural Science Foundation of China (12101011). W. Zhu is supported by the National Natural Science Foundation of China (12201118) and Guangdong
Basic and Applied Basic Research Foundation (2021A1515111018).

\section*{Conflict of interest}
The authors declare that they have no conflict of interest.

\section*{Data Availability} Data sharing is not applicable to this article as no new data were created or analyzed in this study.

\end{document}